**Probability distribution for the Gaussian curvature of the zero level surface of a random function.**

J. H. Hannay,
H. H. Wills Physics Laboratory, University of Bristol,
Tyndall Avenue, Bristol BS8 1TL.  UK.

*Abstract*
A rather natural construction for a smooth random surface in space is the level surface of value zero, or 'nodal' surface $f(x,y,z)=0$, of a (real) random function $f$; the interface between positive and negative regions of the function.  A physically significant local attribute at a point of a curved surface is its Gaussian curvature (the product of its principal curvatures) because, when integrated over the surface it gives the Euler characteristic.  Here the probability distribution for the Gaussian curvature at a random point on the nodal surface $f=0$ is calculated for a statistically homogeneous ('stationary') and isotropic zero mean Gaussian random function $f$.  Capitalizing on the isotropy, a 'fixer' device for axes supplies the probability distribution directly as a multiple integral.  Its evaluation yields an explicit algebraic function with a simple average.  Indeed, this average Gaussian curvature has long been known.  For a non-zero level surface instead of the nodal one, the probability distribution is not fully tractable, but is supplied as an integral expression.

1. *Introduction*
A rather natural construction for a smooth random surface throughout space is the level surface of value zero, or 'nodal' surface $f(x,y,z)=0$, of a statistically homogeneous ('stationary'), isotropic, real, zero mean, Gaussian random function $f$.  Provided the correlation function of $f$ is smooth, the surface is generically smoothly curved and non-self-intersecting (Fig 1), and the regions $f>0$ and $f<0$ on either side of the surface are statistically identical.   It might characterize an interface between two equally populated media or phases in condensed matter, or the neutrally dense surface of density fluctuations in the early universe [Gott et al 1986].  Of interest physically may be the surface's topology, and the most relevant local attribute is then the Gaussian curvature (the product of its two principal curvatures), which varies over the surface.  The integral of the Gaussian curvature over a surface is its Euler characteristic.

The random function $f$ just described is a superposition of an infinity ($\infty$) of plane waves, with random directions and phases, (and amplitudes $\sim 1/\sqrt{\infty}$), fully characterized by its power spectrum, the Fourier transform of its correlation function. (In a condensed matter context the general level surfaces $f=constant$ are sometimes called 'Cahn's model' for random surfaces [Cahn 1965]).

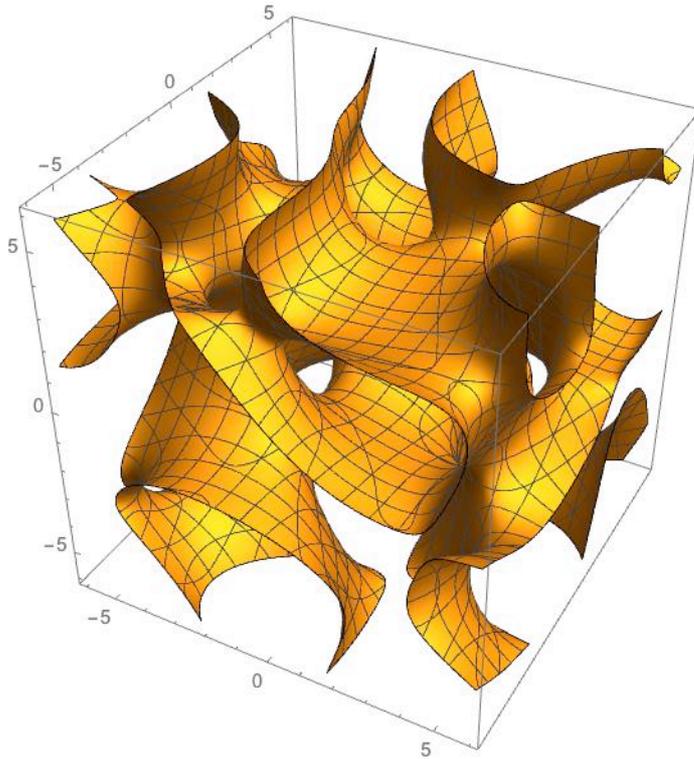

Fig 1. Nodal surface *f*=0 of a random real function *f(x,y,z)* (the portion of it inside a box is shown). The function *f* is statistically homogeneous and isotropic. For the case chosen the constituent waves have equal wavelength (sphere, or hollow ball, power spectrum) and equal amplitudes, and uniformly random directions and phases. Properly there should be infinity of them, here approximated by 50. The convoluted appearance, termed spongelike [Gott et al 1986], shows the predominance of negative Gaussian curvature, saddle-like, surface.

Here the probability distribution for the Gaussian curvature at a (uniformly) random point on the nodal surface (*f*=0) of such a random function *f* will be calculated. The probability distribution turns out to be an elementary (algebraic) function, depending on the two parameters $\langle (\nabla f)^2 \rangle / \langle f^2 \rangle$ and $\langle (\nabla^2 f)^2 \rangle / \langle f^2 \rangle$. These two parameters related respectively to the second and fourth moments of the power spectrum. If the width of the distribution is imagined standardised, a single parameter suffices. These same two basic parameters also appear in other statistical calculations for wavefields, for example [Longuet-Higgins 1958], [Berry and Dennis 2000], [Hannay 2017]. The result is easily generalized (Appendix) to supply the probability distribution for the Gaussian curvature of non-zero level surfaces *f=constant*, but the expression is not fully tractable; one integral remains.

A long known result [Adler and Hasofer 1976] [Gott et al 1986] [Teubner 1991] [Adler and Taylor 2007] in this circumstance supplies average Gaussian curvature, not only for the nodal surface *f*=0 but for any other chosen level surface value *f=constant*. (The terminology is sometimes different: the 'mean Euler characteristic per unit volume' is simply the average Gaussian curvature times the average surface area per unit volume). Strikingly, for the nodal surface, the average is negative, giving a heavily tunnelled, appearance (Fig1).

This negative average can be made plausible algebraically as follows. Let the two principal curvatures of the nodal surface, (with some chosen sign convention, say positive for curved towards higher *f*), be denoted *p* and *q* with *p*>*q*. Then by the ± symmetry of the function *f*, the average <*p*+*q*> is zero. Thus <*p*> = −<*q*>, (both sides are positive, rather than zero, since *p*>*q* by definition). Thus <*p*><*q*> < 0, which makes it reasonable (though is not a proof, of course) that <*pq*> < 0 also. It is usually easier to calculate averages rather than complete probability distributions – and this has been the goal of the previous calculational recipes by others. The direct calculation of this, long known, average Gaussian curvature, by the method here, is indeed a shorter calculation than obtaining the average indirectly from the probability distribution, and is presented in the appendix.

Another known result [Teubner 1991] supplies the probability distribution for the other important type of curvature, the 'mean' curvature *H* (half the sum of the principle curvatures). For the nodal surface *f*=0 the distribution is a very simple expression ($\propto (const + H^2)^{-5/2}$), a symmetric function, as it must be by the ± statistical symmetry of the random function *f*. For level surfaces of non-zero value *f*=*constant*, Teubner obtains the distribution of *H* as an integral expression. The mean and mean square values of *H* are calculated explicitly.

As a final historical note, perhaps the first calculation of the Gaussian curvature probability distribution of a random surface was for a nearly flat landscape function whose indefinitely small height is a random function over the plane [Longuet-Higgins 1958]. The result for the distribution is an integral expression, but the moments can be evaluated. In particular, in contrast to the present nodal surface curvature distribution, the average Gaussian curvature is zero, as it must be, since the landscape has Euler characteristic zero, rather than negative.

The statistical homogeneity and isotropy of *f* imply relations between the averages of pair products of partial derivatives. If indicated by subscripts, these derivatives (evaluated at the same point) average to zero if any of the following apply: *x* subscript total is odd; same for *y*; same for *z*. For example <*ff*$_x$> is zero but <*f*$_{xx}$*f*$_{yy}$> is not. Also, any individual subscript in a pair product average can be shifted to the other member of the pair with the acquisition of a minus sign. For example <*ff*$_{xx}$>=−<*f*$_x^2$> and, (applied twice) <*f*$_{xx}$*f*$_{yy}$>=<*f*$_{xy}^2$>. Zero pair average implies the independence of the members of the pair (at their common point of evaluation), thus *f* and *f*$_x$ are independent.

Capitalizing on the isotropy, a 'fixer' device for axes [Hannay 2017] is to be used to facilitate the calculation. The 'fixer' is the identity $\langle \bullet \rangle = \langle \bullet 2\pi \delta(f_x / f_z) \delta(f_y / f_z) \rangle$, where the dot stands for any intrinsic physical quantity (that has meaning independent of coordinate frame orientation). This device supplies the Gaussian curvature probability distribution directly as a multiple integral; from then on it is a fairly routine matter of evaluation. The result is an explicit algebraic expression for the probability distribution, with simple averages. For the present application the reasoning behind the 'fixer' runs, in summary, as follows. For a location on the nodal surface *f*=0 where the tangent plane happens to be perpendicular to the *z* direction, the Gaussian curvature *K* is $\left( f_{xx} f_{yy} - f_{xy}^2 \right) / f_z^2$.

The requirement that indeed that the nodal surface passes through the chosen point is to be imposed by $\delta(f)|f_z|$ (the 'Jacobian' factor here makes the normal integral of this combination across the surface unity). By statistical isotropy, the z direction is typical of all directions. The $2\pi$ in the identity above is the total solid angle of directions divided by 2 to account for the two possible 'ways up' of the plane.

## 2. The probability distribution of Gaussian curvature

The probability distribution for the Gaussian curvature $K$ (the product of the two principal curvatures) at a random point on the nodal surface $f=0$, is given as follows. The angle bracket averages refer to ensemble averages, or equivalently, by homogeneity, volume averages. The denominator is the average nodal surface area per unit volume. The integrals are along the entire real axis unless otherwise specified.

$$P(K) = \frac{2\pi\left\langle \delta(f)\,|f_z|\,\delta\!\left(\frac{f_x}{f_z}\right)\delta\!\left(\frac{f_y}{f_z}\right)\delta\!\left(K - \frac{f_{xx}f_{yy}-f_{xy}^2}{f_z^2}\right)\right\rangle}{2\pi\left\langle \delta(f)\,|f_z|\,\delta\!\left(\frac{f_x}{f_z}\right)\delta\!\left(\frac{f_y}{f_z}\right)\right\rangle} \qquad (1)$$

$$= \frac{2\pi\langle\delta(f_x)\rangle\langle\delta(f_y)\rangle\left\langle \delta(f)\,|f_z|^5\,\delta\!\left(Kf_z^2 - (f_{xx}f_{yy}-f_{xy}^2)\right)\right\rangle}{2\pi\langle|f_z|^3\rangle\langle\delta(f_x)\rangle\langle\delta(f_y)\rangle\langle\delta(f)\rangle} \qquad (2)$$

$$= \frac{\left\langle \delta(f)\,|f_z|^5\,\delta\!\left(Kf_z^2 - (f_{xx}f_{yy}-f_{xy}^2)\right)\right\rangle}{\langle|f_z|^3\rangle\langle\delta(f)\rangle} \qquad (3)$$

$$= \left[\left(\tfrac{1}{2}\Gamma_{z,z}\right)^2 \bigg/ \sqrt{\frac{\Gamma_{z,z}}{2\pi}}\sqrt{\frac{\Gamma_{0,0}}{2\pi}}\right]\left\langle \delta(f)\,|f_z|^5 \int \exp(i\gamma Kf_z^2)\exp(-i\gamma(f_{xx}f_{yy}-f_{xy}^2))\frac{d\gamma}{2\pi}\right\rangle \qquad (4)$$

Here, and in the next steps, various definitions have been, or will be, used:

$$\Gamma_{0,0} = \frac{1}{\langle f^2\rangle},\quad \Gamma_{z,z} = \frac{1}{\langle f_z^2\rangle}(=\Gamma_{x,x}=\Gamma_{y,y}) = -\frac{1}{\langle ff_{xx}\rangle} = -\frac{1}{\langle ff_{yy}\rangle},\quad \Gamma_{xy,xy} = \frac{1}{\langle f_{xy}^2\rangle} = \frac{1}{\langle f_{xx}f_{yy}\rangle}. \qquad (5)$$

Also required are a 3x3 matrix $\boldsymbol{\Gamma}$ (the *inverse* of a matrix of pair product averages) and its 2x2 lower right submatrix $\boldsymbol{\Gamma}'$, whose inverse is therefore a more complicated expression in terms of product averages:

$$\Gamma = \begin{pmatrix} \langle f^2 \rangle & \langle ff_{xx} \rangle & \langle ff_{yy} \rangle \\ \langle ff_{xx} \rangle & \langle f_{xx}^2 \rangle & \langle f_{xx}f_{yy} \rangle \\ \langle ff_{yy} \rangle & \langle f_{xx}f_{yy} \rangle & \langle f_{yy}^2 \rangle \end{pmatrix}^{-1} \tag{6}$$

$$\Gamma' \equiv \begin{pmatrix} \Gamma_{xx,xx} & \Gamma_{xx,yy} \\ \Gamma_{xx,yy} & \Gamma_{yy,yy} \end{pmatrix} = \begin{pmatrix} \langle f_{xx}^2 \rangle - \frac{\langle ff_{xx} \rangle^2}{\langle f^2 \rangle} & \langle f_{xx}f_{yy} \rangle - \frac{\langle ff_{xx} \rangle \langle ff_{yy} \rangle}{\langle f^2 \rangle} \\ \langle f_{xx}f_{yy} \rangle - \frac{\langle ff_{xx} \rangle \langle ff_{yy} \rangle}{\langle f^2 \rangle} & \langle f_{yy}^2 \rangle - \frac{\langle ff_{yy} \rangle^2}{\langle f^2 \rangle} \end{pmatrix}^{-1} \tag{7}$$

Resuming the calculation:

$$P(K) = \frac{\left(\frac{1}{2}\Gamma_{z,z}\right)^2}{\sqrt{\frac{\Gamma_{z,z}}{2\pi}}\sqrt{\frac{\Gamma_{0,0}}{2\pi}}} \int \langle |f_z|^5 \exp(i\gamma K f_z^2) \rangle \langle \exp(i\gamma(f_{xy}^2)) \rangle \langle \delta(f) \exp(-i\gamma f_{xx}f_{yy}) \rangle \frac{d\gamma}{2\pi} \tag{8}$$

$$= \left[\left(\tfrac{1}{2}\Gamma_{z,z}\right)^2 \Big/ \sqrt{\tfrac{\Gamma_{z,z}}{2\pi}}\sqrt{\tfrac{\Gamma_{0,0}}{2\pi}}\right] \int \frac{d\gamma}{2\pi} \int_0^\infty (f_z^2)^2 \exp(i\gamma K f_z^2)\exp\left(-\tfrac{1}{2}\Gamma_{z,z}(f_z^2)\right) d(f_z^2) \times$$

$$\times \sqrt{\tfrac{\Gamma_{xy,xy}}{2\pi}} \int \exp(i\gamma f_{xy}^2)\exp\left(-\tfrac{1}{2}\Gamma_{xy,xy}f_{xy}^2\right) df_{xy} \times \tag{9}$$

$$\times \sqrt{\tfrac{\det \Gamma}{(2\pi)^3}} \iint \exp(-i\gamma f_{xx}f_{yy})\exp\left[-\tfrac{1}{2}(0,f_{xx},f_{yy})\Gamma\begin{pmatrix}0\\f_{xx}\\f_{yy}\end{pmatrix}\right] df_{xx}\, df_{yy}$$

$$= \left[\left(\tfrac{1}{2}\Gamma_{z,z}\right)^2 \Big/ \sqrt{\tfrac{\Gamma_{z,z}}{2\pi}}\sqrt{\tfrac{\Gamma_{0,0}}{2\pi}}\right] \int \frac{d\gamma}{2\pi}\sqrt{\tfrac{\Gamma_{z,z}}{2\pi}}\frac{2!}{\left(\tfrac{1}{2}\Gamma_{z,z}-i\gamma K\right)^3} \times \frac{1}{\sqrt{1-2i\gamma/\Gamma_{xy,xy}}} \times$$

$$\times \sqrt{\tfrac{\det \Gamma}{(2\pi)^3}} \iint \exp(-i\gamma f_{xx}f_{yy})\exp\left[-\tfrac{1}{2}(f_{xx},f_{yy})\Gamma'\begin{pmatrix}f_{xx}\\f_{yy}\end{pmatrix}\right] df_{xx}\, df_{yy} \tag{10}$$

$$= \left(\left(\tfrac{1}{2}\Gamma_{z,z}\right)^2 \Big/ \sqrt{\tfrac{\Gamma_{0,0}}{2\pi}}\right)\sqrt{\tfrac{\det \Gamma}{(2\pi)^3}} \times$$

$$\times \int \frac{d\gamma}{2\pi}\frac{2!}{\left(\tfrac{1}{2}\Gamma_{z,z}-i\gamma K\right)^3} \times \frac{1}{\sqrt{1-2i\gamma/\Gamma_{xy,xy}}} \times \sqrt{\frac{(2\pi)^2}{\Gamma_{xx,xx}\Gamma_{yy,yy}-(\Gamma_{xx,yy}+i\gamma)^2}} \tag{11}$$

$$= \left(\left(\tfrac{1}{2}\Gamma_{z,z}\right)^2 \middle/ \sqrt{\frac{\Gamma_{0,0}}{2\pi}}\right) \sqrt{\frac{\det \Gamma}{(2\pi)^3}} \int \frac{d\gamma}{2\pi} \frac{2!}{\left(\tfrac{1}{2}\Gamma_{z,z} - i\gamma K\right)^3} \times$$

$$\times \frac{1}{\sqrt{1 - 2i\gamma/\Gamma_{xy,xy}}} \times \sqrt{\frac{2\pi}{\Gamma_{xx,xx} - \Gamma_{xx,yy} - i\gamma}} \times \sqrt{\frac{2\pi}{\Gamma_{xx,xx} + \Gamma_{xx,yy} + i\gamma}} \tag{12}$$

This would be expected to be an Elliptic integral, but remarkably the first two square roots in the integrand are proportional to each other - their branch points coincide - by virtue of a relation (to be proved below):

$$\Gamma_{xx,xx} - \Gamma_{xx,yy} = \tfrac{1}{2}\Gamma_{xy,xy}. \tag{13}$$

The consequence is that the integrand has only the final square root remaining, its other singularities being poles, and is therefore straightforward to evaluate by residues. (This same coincidence of singularities similarly simplified a vortex reconnection rate calculation for a random complex function [Hannay 2017]).

$$P(K) = \left\{ \left(\left(\tfrac{1}{2}\Gamma_{z,z}\right)^2 \middle/ \sqrt{\frac{\Gamma_{0,0}}{2\pi}}\right) \sqrt{\frac{\det \Gamma}{(2\pi)^3}} \sqrt{\tfrac{1}{2}\Gamma_{xy,xy}}\, 2! \right\} \times$$

$$\int \frac{1}{\left(\tfrac{1}{2}\Gamma_{z,z} - i\gamma K\right)^3} \frac{1}{\Gamma_{xx,xx} - \Gamma_{xx,yy} - i\gamma} \frac{1}{\sqrt{\Gamma_{xx,xx} + \Gamma_{xx,yy} + i\gamma}}\, d\gamma \tag{14}$$

Before quoting the result, the proof of the relation (13) is as follows. The isotropy of the power spectrum means that, using standard spherical polar coordinates,

$$\langle f_x^2 \rangle = \langle -ff_{xx} \rangle = \langle f^2 \rangle \overline{k^2} \tfrac{1}{4\pi} \int_0^{2\pi} \int_{-1}^{1} \cos^2\phi \cos^2\theta \sin\theta\, d\theta\, d\phi = \tfrac{1}{3}\langle f^2 \rangle \overline{k^2}$$

$$\langle f_{xx}^2 \rangle = \langle ff_{xxxx} \rangle = \langle f^2 \rangle \overline{k^4} \tfrac{1}{4\pi} \int_0^{2\pi} \int_{-1}^{1} \cos^4\phi \cos^4\theta \sin\theta\, d\theta\, d\phi = \tfrac{1}{5}\langle f^2 \rangle \overline{k^4} \text{ and} \tag{15}$$

$$\langle f_{xx} f_{yy} \rangle = \langle f_{xy}^2 \rangle = \langle ff_{xxyy} \rangle = \langle f^2 \rangle \overline{k^4} \tfrac{1}{4\pi} \int_0^{2\pi} \int_{-1}^{1} \cos^2\phi \sin^2\phi \cos^4\theta \sin\theta\, d\theta\, d\phi = \tfrac{1}{15}\langle f^2 \rangle \overline{k^4},$$

where $k$ is wavenumber in the power spectrum, and the straight overbar means average over the power spectrum. The $\Gamma'^{-1}$ matrix is therefore

$$(\Gamma')^{-1} = \begin{pmatrix} \Gamma_{xx,xx} & \Gamma_{xx,yy} \\ \Gamma_{xx,yy} & \Gamma_{xx,xx} \end{pmatrix}^{-1} = \frac{\langle f^2 \rangle \overline{k^4}}{15} \begin{pmatrix} 3 - 5\rho/3 & 1 - 5\rho/3 \\ 1 - 5\rho/3 & 3 - 5\rho/3 \end{pmatrix} \tag{16}$$

where $\rho \equiv \overline{k^2}^2 / \overline{k^4} / \leq 1$. The difference of the two elements of the top row of the inverse of this matrix, $\Gamma_{xx,xx} - \Gamma_{xx,yy}$, is $\frac{1}{2} / \left(\frac{1}{15}\langle f^2 \rangle \overline{k^4}\right)$ (independent of $\rho$) which equals $\frac{1}{2} / \langle f_{xy}^2 \rangle = \frac{1}{2} \Gamma_{xy,xy}$, as required.

Returning to the evaluation of P(K), the branch point is on the negative imaginary axis, and the simple pole is on the positive imaginary axis. By deforming the contour of integration upwards the residue contribution from this simple pole is acquired. The triple pole may or may not contribute according to the sign of K. If K is positive then it does, since its location is on the positive imaginary axis, if K is negative, its location has switched sign and therefore it does not contribute. Denoting $a \equiv \frac{1}{2}\Gamma_{z,z}/\Gamma_{0,0}$, $b \equiv (\Gamma_{xx,xx} - \Gamma_{xx,yy})/\Gamma_{0,0}$, $c \equiv (\Gamma_{xx,xx} + \Gamma_{xx,yy})/\Gamma_{0,0}$ (all of which are positive), the prefactor in the brace in (14) becomes $\{a^2 b \sqrt{c}/\pi\}$ and

$$P(K) = \left\{\frac{a^2 b \sqrt{c}}{\pi}\right\} \int \frac{1}{(a-i\gamma K)^3} \frac{1}{b-i\gamma} \frac{1}{\sqrt{c+i\gamma}} d\gamma \quad = \frac{2\pi}{(a-bK)^3 \sqrt{b+c}} \quad \text{if } K<0,$$

or (17)

$$= \frac{\pi\sqrt{K}\left(-15a^2 + 10a(b-2c)K + (-3b^2 + 4bc - 8c^2)K^2\right)}{4(a-Kb)^3(a+Kc)^{5/2}} + \frac{2\pi}{(a-Kb)^3\sqrt{b+c}} \quad \text{if } K>0.$$

(18)

(In this latter case K>0, there are singularities in the two terms at K=a/b but they cancel).

These formulas, (17) and (18), are the result. With the definitions above for a, b, c, and the substitutions from (5) and (16), namely

$$\Gamma_{z,z} = 3 / \left(\langle f^2\rangle \sqrt{\overline{k^4}\rho}\right), \quad \Gamma_{xx,xx} = \left(\overline{k^4}\langle f^2\rangle\right)^{-1} \frac{15}{4}\left(\frac{9-5\rho}{6-5\rho}\right), \quad \Gamma_{xx,yy} = \left(\overline{k^4}\langle f^2\rangle\right)^{-1} \frac{15}{4}\left(\frac{-3+5\rho}{6-5\rho}\right)$$

(19)

(recalling $\rho \equiv \overline{k^2}^2 / \overline{k^4}$), the probability distribution P(K) is expressed in terms of the power spectrum. This is graphed in fig 2. The average value of K, $\int_{-\infty}^{\infty} K\, P(K)\, dK$ is long known, $-\overline{k^2}/6$ [Adler and Hasofer 1976]. Since P(K) decays as $|K|^{-3}$ for large K, the mean square of K diverges logarithmically.

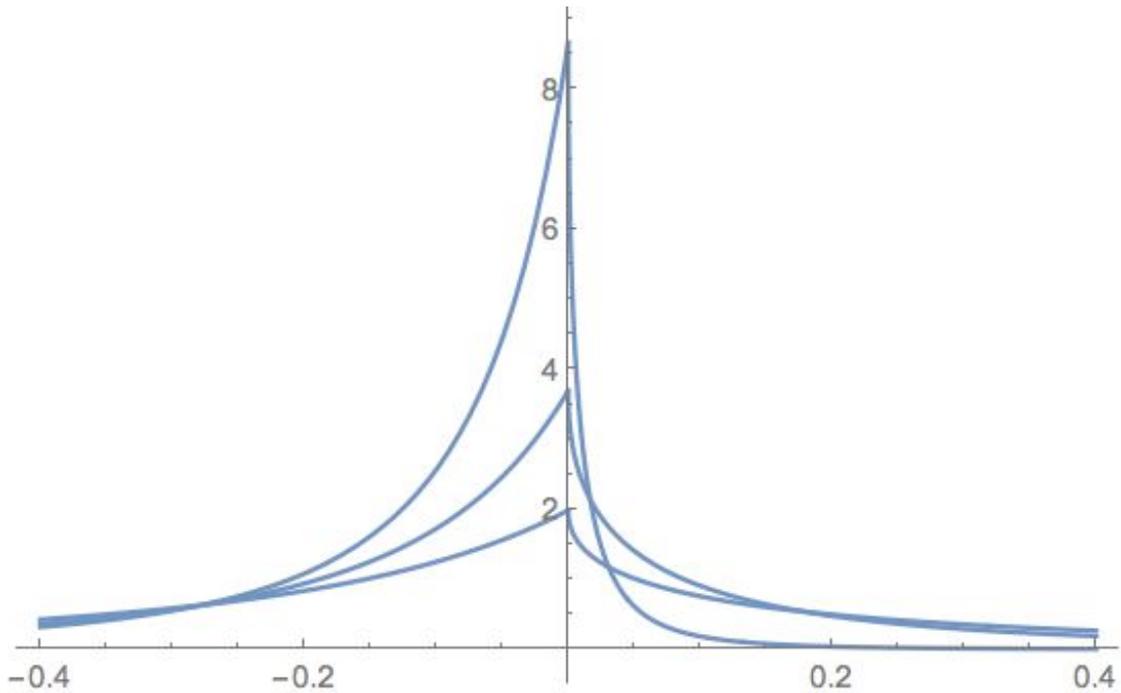

Fig 2. The probability distribution $P(K)$ of the Gaussian curvature at a random position on the nodal surface $f=0$ of an homogeneous isotropic random function. The distribution is shown for the three values 1 (tallest graph, height $5\sqrt{3}$), 1/3, 1/9 of the dimensionless parameter, $\rho \equiv \overline{k^2}^2 / \overline{k^4} \leq 1$, where the overbar means average over the power spectrum. The scale width for the graph is fixed by taking $\overline{k^4} = 1$. The graphs are noticeably asymmetric and the average value of $K$ is $-\overline{k^2}/6 \ (= -\sqrt{\rho}/6)$, a long known result.

*Appendix. Generalizing to a non-zero level surface f=F*

The probability distribution $P_F(K)$ for the Gaussian curvature at a random point on the surface $f=F$ is calculated similarly. The equations with notable changes are as follows. An extra Gaussian factor in $F$ appears in (A2) and then cancels again in (A4), and in (A3) the previous zeros in the vector components have been replaced by $F$. An extra $F$ dependent exponential in (A4) leads to an extra one in (A5). The result (A6) is simply (14) but with an extra exponential factor that prevents its analytical evaluation.

$$P_F(K) = \frac{\langle \delta(f-F) \ |f_z|^5 \ \delta(Kf_z^2 - (f_{xx}f_{yy} - f_{xy}^2))\rangle}{\langle |f_z|^3 \rangle \langle \delta(f-F)\rangle} \tag{A1}$$

$$= \left[(\tfrac{1}{2}\Gamma_{z,z})^2 \Big/ \sqrt{\frac{\Gamma_{z,z}}{2\pi}}\sqrt{\frac{\Gamma_{0,0}}{2\pi}}\exp(-\tfrac{1}{2}\Gamma_{0,0}F^2)\right] \times$$
$$\left\langle \delta(f-F) \ |f_z|^5 \int \exp(i\gamma Kf_z^2)\exp(-i\gamma(f_{xx}f_{yy}-f_{xy}^2))\frac{d\gamma}{2\pi}\right\rangle \tag{A2}$$

$$= \left[(\tfrac{1}{2}\Gamma_{z,z})^2 \Big/ \sqrt{\frac{\Gamma_{z,z}}{2\pi}}\sqrt{\frac{\Gamma_{0,0}}{2\pi}}\exp(-\tfrac{1}{2}\Gamma_0 F^2)\right]\times$$
$$\times\int \frac{d\gamma}{2\pi}\int_0^\infty (f_z^2)^2 \exp(i\gamma Kf_z^2)\exp(-\tfrac{1}{2}\Gamma_{z,z}(f_z^2))d(f_z^2)\times$$
$$\times\sqrt{\frac{\Gamma_{xy,xy}}{2\pi}}\int \exp(i\gamma f_{xy}^2)\exp(-\tfrac{1}{2}\Gamma_{xy,xy}f_{xy}^2)df_{xy}\times \tag{A3}$$
$$\times\sqrt{\frac{\det\Gamma}{(2\pi)^3}}\iint \exp(-i\gamma f_{xx}f_{yy})\exp\left[-\tfrac{1}{2}(F,f_{xx},f_{yy})\Gamma\begin{pmatrix}F\\f_{xx}\\f_{yy}\end{pmatrix}\right]df_{xx}\,df_{yy}$$

$$= \left[(\tfrac{1}{2}\Gamma_{z,z})^2 \Big/ \sqrt{\frac{\Gamma_{z,z}}{2\pi}}\sqrt{\frac{\Gamma_{0,0}}{2\pi}}\right]\times \int\frac{d\gamma}{2\pi}\int_0^\infty (f_z^2)^2\exp(i\gamma Kf_z^2)\exp(-\tfrac{1}{2}\Gamma_{z,z}(f_z^2))d(f_z^2)\times$$
$$\times\sqrt{\frac{\Gamma_{xy,xy}}{2\pi}}\int\exp(i\gamma f_{xy}^2)\exp(-\tfrac{1}{2}\Gamma_{xy,xy}f_{xy}^2)df_{xy}\times\sqrt{\frac{\det\Gamma}{(2\pi)^3}}\times \tag{A4}$$
$$\times\iint\exp(-i\gamma f_{xx}f_{yy})\exp(-F\Gamma_{0,xx}(f_{xx},f_{yy}).(1,1))\,\exp\left[-\tfrac{1}{2}(f_{xx},f_{yy})\Gamma'\begin{pmatrix}f_{xx}\\f_{yy}\end{pmatrix}\right]df_{xx}\,df_{yy}$$

$$= \left( \left( \tfrac{1}{2}\Gamma_{z,z}\right)^2 \Big/ \sqrt{\frac{\Gamma_{0,0}}{2\pi}}\right) \sqrt{\frac{\det \Gamma}{(2\pi)^3}} \times$$

$$\times \int \frac{d\gamma}{2\pi} \frac{2!}{\left(\tfrac{1}{2}\Gamma_{z,z} - i\gamma K\right)^3} \times \frac{1}{\sqrt{1 - 2i\gamma/\Gamma_{xy,xy}}} \times \sqrt{\frac{(2\pi)^2}{\Gamma_{xx,xx}\Gamma_{yy,yy} - (\Gamma_{xx,yy} + i\gamma)^2}} \quad (A5)$$

$$\times \exp\left[ -\frac{1}{2}(F\Gamma_{0,xx})^2 (1,1) \begin{pmatrix} \Gamma_{xx,xx} & \Gamma_{xx,yy} + i\gamma \\ \Gamma_{xx,yy} + i\gamma & \Gamma_{xx,xx} \end{pmatrix}^{-1} \begin{pmatrix} 1 \\ 1 \end{pmatrix} \right]$$

$$= \left( \left( \tfrac{1}{2}\Gamma_{z,z}\right)^2 \Big/ \sqrt{\frac{\Gamma_{0,0}}{2\pi}}\right) \sqrt{\frac{\det \Gamma}{(2\pi)^3}} 2! \sqrt{\tfrac{1}{2}\Gamma_{xy,xy}}$$

$$\times \int d\gamma \frac{1}{\left(\tfrac{1}{2}\Gamma_{z,z} - i\gamma K\right)^3} \frac{1}{\Gamma_{xx,xx} - \Gamma_{xx,yy} - i\gamma} \frac{1}{\sqrt{\Gamma_{xx,xx} + \Gamma_{xx,yy} + i\gamma}} \times \quad (A6)$$

$$\times \exp\left[ -(F\Gamma_{0,xx})^2 \big/ (\Gamma_{xx,xx} + \Gamma_{xx,yy} + i\gamma) \right]$$

where $\Gamma_{0,xx} = \langle f_x^2 \rangle \Big/ \left( \langle f^2 \rangle \langle f_{xx}^2 \rangle + \langle f^2 \rangle \langle f_{xx} f_{yy} \rangle - 2\langle f_x^2 \rangle^2 \right) = \overline{k^2} \Big/ \left( \langle f^2 \rangle \overline{k^4}(\tfrac{4}{3} - 2\rho) \right)$ and the other constants are supplied, as before, by (5), (15) and (19). The average value of $K$, $\int_{-\infty}^{\infty} K\, P_F(K)\, dK$, is again long known $\tfrac{1}{6}\overline{k^2}\left(F^2/\langle f^2 \rangle - 1\right)$ [Adler and Hasofer 1976]. Here it would be most easily accessed by direct calculation; by replacing the numerator of (A1) by $\langle \delta(f - F)\ |f_z|(f_{xx}f_{yy} - f_{xy}^2) \rangle$ and proceeding.